\input amstex
\magnification=1200
\documentstyle{amsppt}
\NoRunningHeads
\NoBlackBoxes
\topmatter
\title Tactical games\linebreak \& \linebreak behavioral self-organization
\endtitle
\author Denis V. Juriev\endauthor
\affil ul.Miklukho-Maklaya 20-180, Moscow 117437 Russia\linebreak
(e-mail: denis\@juriev.msk.ru)\endaffil
\date adap-org/9911147\newline November, 15, 1999\enddate
\abstract\nofrills The interactive game theoretical approach to tactics and
behavioral self-organization is developed. Though it uses the interactive 
game theoretical formalization of dialogues as psycholinguistic phenomena, 
the crucial role is played by the essentially new concept of a tactical game.
Applications to the perception processes and related subjects
(memory, recollection, image understanding, imagination) are discussed
together with relations to the computer vision and pattern recognition (the
dynamical formation of patterns and perception models during perception as
a result of its self-organization) and computer games (modelling of the
tactical behavior and self-organization, tactical RPG and elaboration of new
tactical game techniques). The appendix is devoted to the operative computer
games and the user programming of operative units in a multi-user online
operative computer game.
\endabstract
\toc
\specialhead Introduction\endspecialhead
\head I. Interactive games\endhead
\subhead 1.1. Interactive systems and games\endsubhead
\subhead 1.2. The $\varepsilon$-representation\endsubhead
\subhead 1.3. Hidden interactivity of ordinary differential games\endsubhead
\head II. Dialogues and verbalizable games\endhead
\subhead 2.1. Dialogues\endsubhead
\subhead 2.2. The verbalization of interactive games\endsubhead
\subhead 2.3. Psycholinguistic encoding/decoding of dialogues\endsubhead
\head III. \ Tactical games\endhead
\subhead 3.1. Parametric interactive games and external controls\endsubhead
\subhead 3.2. Dialogue and comments\endsubhead
\subhead 3.3. Tactical games\endsubhead
\subhead 3.4. Applications to computer games (tactical action games and the
modelling of tactical behavior)\endsubhead
\head IV. \ Behavioral self-organization\endhead
\subhead 4.1. Tactical behavioral self-organization\endsubhead
\subhead 4.2. Applications to computer games (artificial intelligence and
tactical self-organization, tactical RPG)\endsubhead
\head V. Tactics and behavioral self-organization for perception
games\endhead
\subhead 5.1. Perception games and kaleidoscope-roulettes\endsubhead
\subhead 5.2. Tactics and behavioral self-organization for perception
games. Perceptive oracles\endsubhead
\head VI. \ Tactical interaction, tactical synthesis and tactical localization
of control systems\endhead
\subhead 6.1. Tactical interaction of control systems\endsubhead
\subhead 6.2. Tactical synthesis of control systems\endsubhead
\subhead 6.3. Tactical localization of control systems\endsubhead
\head Appendix \ \ \ \ \ \ \ \ \ \ A. Operative computer games\endhead
\subhead A.1. Tactics and operative art. Operative computer games\endsubhead
\subhead A.2. Programming of operative units in operative computer
games\endsubhead
\endtoc
\endtopmatter
\document

\specialhead Introduction\endspecialhead

The mathematical formalism of interactive games, which extends one of 
ordinary games (see e.g.[1]) and is based on the concept of an interactive 
control, was recently proposed by the author [2] to take into account 
the complex composition of controls of a real human person, which are often 
complicated couplings of his/her cognitive and known controls with 
the unknown subconscious behavioral reactions. The goal of this article 
is to describe the tactical phenomena and behavioral self-organization 
in terms of interactive games. Though the crucial role is played by 
the interactive game theoretical concepts of dialogues 
and verbalizable games [3], the proposed constructions are essentially new
as ideologically as technically. However, it should be specially emphasized
that tactics is derived from the interactivity, and roughly speaking, the 
first may be thought as just an art to manipulate the Unknown, which manifests 
itself in the least, without making it known completely\footnote{\ To my mind 
this is very reasonable if the ecological approach to the Unknown is adopted 
and we are not willing the least to be written by the golden letters of our 
Knowledge into the Red Book of the Universe.\newline}. Otherwords, besides the
reason tactics includes the sphere of ``another effort of mind'', which is 
important not only in practice but also in theory\footnote{\ The simplest 
explanation why the transition from interactivity to tactics is necessary is 
that the description of interactive process is a tactical procedure (see 
e.g.[4]), otherwords, interactive games have their descriptors among tactical 
games whereas the class of tactical games contains all their descriptors 
and therefore is self-consistent. It is essential to mark that {\sl the 
description of any interactive phenomena is interactive itself}.\newline}. 
And because this sphere is important for all moments of our life it is not 
strange that virtually all known forms of human activity such as scientific
researches and economics, sport or military actions, fine or martial arts, 
literature, music and theatre, psychotherapy and even magic may be regarded 
as certain tactical games, which therefore have the universal meaning.

It is well-known that the world around us is much deeper than what we see,
are able to perceive or even to comprehend, so it is essential not only to
understand the known forms of human (individual and collective) activity
in terms of tactical games but also to try to unravel the abstract mathematical
foundations of all tactical phenomena and thus to extend the sphere of 
tactical activity presumably into the abstract intelligible direction
(certainly, completely conserving its practical orientation). Thereas the 
first purpose is the goal of the present article, the least will be discussed 
in the forthcoming one devoted to the game theoretical description of 
{\sl dialectics}.

The article is organized in the following manner. Two first paragraphs are
introductory. They are devoted to the general interactive games and the
dialogues, respectively. The third paragraph is devoted to the new concept
of tactical games. Interactive game theoretical models of behavioral
self-organizations are discussed in the fourth paragraph. Other paragraphs
contains special topics: in the firth paragraph the attention is concentrated 
on tactics and behavioral self-organization constructed for the perception 
games and the kaleidoscope-roulettes in particular, the sixth one is devoted
to the tactical interaction, tactical synthesis and tactical localization
of control systems.

The article includes an appendix, where the operative computer games and the 
user programming of operative units in a multi-user online operative computer 
game are discussed. 

The article is organized as series of definitions supplied by the clarifying
remarks, which together describe the entire construction.

\head I. Interactive games\endhead

First of all, let us expose the basic principles of interactive games
in general.

\subhead 1.1. Interactive systems and games\endsubhead

\definition{Definition 1 [2]} An {\it interactive system\/} (with $n$
{\it interactive controls\/}) is a control system with $n$ independent 
controls coupled with unknown or incompletely known feedbacks (the feedbacks
as well as their couplings with controls are of a so complicated nature that 
their can not be described completely). An {\it interactive game\/} is a game 
with interactive controls of each player.
\enddefinition

Below we shall consider only deterministic and differential interactive
systems. In this case the general interactive system may be written in the 
form:
$$\dot\varphi=\Phi(\varphi,u_1,u_2,\ldots,u_n),\tag1$$
where $\varphi$ characterizes the state of the system and $u_i$ are
the interactive controls:
$$u_i(t)=u_i(u_i^\circ(t),\left.[\varphi(\tau)]\right|_{\tau\leqslant t}),$$
i.e. the independent controls $u_i^\circ(t)$ coupled with the feedbacks on
$\left.[\varphi(\tau)]\right|_{\tau\leqslant t}$. One may suppose that the
feedbacks are integrodifferential on $t$.

However, it is reasonable to consider the {\it differential interactive
games}, whose feedbacks are purely differential. It means that
$$u_i(t)=u_i(u_i^\circ(t),\varphi(t),\ldots,\varphi^{(k)}(t)).$$
A reduction of general interactive games to the differential ones via the 
introducing of the so-called {\it intention fields\/} was described in [2]. 
Below we shall consider the differential interactive games only if the opposite 
is not specified explicitely.

The interactive games introduced above may be generalized in the following 
ways. 

The first way, which leads to the {\it indeterminate interactive games},
is based on the idea that the pure controls $u_i^\circ(t)$ and the 
interactive controls $u_i(t)$ should not be obligatory related in the
considered way. More generally one should only postulate that there are
some time-independent quantities $F_\alpha(u_i(t),u_i^\circ(t),\varphi(t),
\ldots,\varphi^{(k)}(t))$ for the independent magnitudes $u_i(t)$ and 
$u_i^\circ(t)$. Such claim is evidently weaker than one of Def.1. For 
instance, one may consider the inverse dependence of the pure and 
interactive controls: $u_i^\circ(t)=u_i^\circ(u_i(t),\varphi(t),\ldots,
\varphi^{(k)}(t))$.

The inverse dependence of the pure and interactive controls has a nice
psychological interpretation. Instead of thinking of our action consisting
of the conscious and unconscious parts and interpreting the least as 
unknown feedbacks which ``dress'' the first, one is able to consider
our action as a single whole whereas the act of consciousness is in
the extraction of a part which it declares as its ``property''. So
interpreted the function $u_i^\circ(u_i,\varphi,\ldots,\varphi^{(k)})$ 
realizes the ``filtering'' procedure.

The second way, which leads to the {\it coalition interactive games}, is
based on the idea to consider the games with coalitions of actions and to
claim that the interactive controls belong to such coalitions. In this case
the evolution equations have the form
$$\dot\varphi=\Phi(\varphi,v_1,\ldots,v_m),\tag2$$
where $v_i$ is the interactive control of the $i$-th coalition. If the 
$i$-th coalition is defined by the subset $I_i$ of all players then
$$v_i=v_i(\varphi(t),\ldots,\varphi^{(k)}(t),u^\circ_j| j\in I_i).$$
Certainly, the intersections of different sets $I_i$ may be non-empty so
that any player may belong to several coalitions of actions. Def.1 gives the
particular case when $I_i=\{i\}$.

The coalition interactive games may be an effective tool for an analysis of
the collective decision making in the real coalition games that spread the
applicability of the elaborating interactive game theory to the diverse 
problems of sociology. 

\remark{Remark 1} One is able to consider interactive games of discrete time 
in the similar manner.
\endremark

\remark{Remark 2} If one suspect that the explicit dependence of the feedbacks
on the derivatives of $\varphi$ is not correct because they are determined
via the evolution equations governed by the interactive controls, it is
reasonable to use the inverse dependence of pure and interactive controls.
\endremark

\subhead 1.2. The $\varepsilon$-representations\endsubhead
Interactive games are games with incomplete information by their nature.
However, this incompleteness is in the unknown feedbacks, not in the 
unknown states. The least situation is quite familiar to specialists in
game theory and there is a lot of methods to have deal with it. For
instance, the unknown states are interpreted as independent controls of
the virtual players and some muppositions on their strategies are done.
To transform interactive games into the games with an incomplete information
on the states one can use the following trick, which is called the 
$\varepsilon$-representation of the interactive game.

\definition{Definition 2} The $\varepsilon$-representation of the 
differential interactive game is a representation of the interactive controls 
$u_i(t)$ in the form
$$u_i(t)=u_i(u^\circ_i(t),\varphi(t),\ldots\varphi^{(k)}(t);\varepsilon_i(t))$$
with the {\sl known\/} function $u_i$ of its arguments $u_i^\circ$,
$\varphi,\ldots,\varphi^{(k)}$ and $\varepsilon_i$, whereas 
$$\varepsilon_i(t)=\varepsilon_i(u^\circ_i(t),\varphi(t),\ldots,\varphi^{(k)}
(t)$$ 
is the {\sl unknown\/} function of $u_i^\circ$ and $\varphi,\ldots,
\varphi^{(k)}$.
\enddefinition

\remark{Remark 3} The derivatives of $\varphi$ may be excluded from the
feedbacks in the way described above. 
\endremark

\remark{Remark 4} One is able to consider the $\varepsilon$-representations
of the indeterminate and coalition interactive games.
\endremark

$\varepsilon_i$ are interpreted as parameters of feedbacks and, thus, 
characterize the internal {\sl existential\/} states of players. It motivates 
the notation $\varepsilon$. Certainly, $\varepsilon$-parameters are not
really states being the unknown functions of the states and pure controls,
however, one may sometimes to apply the standard procedures of the theory
of games with incomplete information on the states. For instance, it is
possible to regard $\varepsilon_i$ as controls of the virtual players.
The na{\"\i}vely introduced virtual players only double the ensemble of 
the real ones in the interactive games but in the coalition interactive 
games the collective virtual players are observed. More sophisticated
procedures generate ensembles of virtual players of diverse structure.

Precisely, if the derivatives of $\varphi$ are excluded from the feedbacks (at 
least, from the interactive controls $u_i$ as functions of the pure controls,
states and the $\varepsilon$-parameters) the evolution equation will have 
the form
$$\dot\varphi(t)=\Phi(\varphi,u_1(u^\circ_1(t),\varphi(t);\varepsilon_1(t)),
\ldots,u_n(u_n^\circ(t),\varphi(t);\varepsilon_n(t))),\tag3$$
so it is consistent to regard the equations as ones of the controlled
system with the ordinary controls $u_1,\ldots,u_n,\varepsilon_1,\ldots,
\varepsilon_n$. One may consider a new game postulating that these
controls are independent. Such game will be called the ordinary differential
game associated with the $\varepsilon$-representation of the interactive
game.

\subhead 1.3. Hidden interactivity of ordinary differential games [5]
\endsubhead Let us consider an arbitrary ordinary differential game with
the evolution equations
$$\dot\varphi=\Phi(\varphi,u_1,u_2,\ldots,u_n),$$
where $\varphi$ characterizes the state of the system and $u_i$ are
the ordinary controls. Let us fix a player. For simplicity of notations we
shall suppose that it is the first one. As a rule the players have their
algorithms of predictions of the behaviour of other players. For a fixed
moment $t_0$ of time let us consider the prediction of the first player for
the game. It consists of the predicted controls $u^\circ_{[t_0];i}(t)$ 
($t>t_0$; $i\ge2$) of all players and the predicted evolution of the system 
$\varphi^\circ_{[t_0]}(t)$. Let us fix $\Delta t$ and consider the real and 
predicted controls for the moment $t_0+\Delta t$. Of course, they may be
different because other players use another algorithms for the game prediction.
One may interpret the real controls $u_i(t)$ ($t=t_0+\Delta t$; $i\ge2$) of 
other players as interactive ones whereas the predicted controls 
$u^\circ_{[t_0];i}(t)$ as pure ones, i.e. to postulate their relation in the 
form:
$$u_i(t)=u_i(u^\circ_{[t-\Delta t];i}(t);\varphi^\circ_{[t_0]}(\tau)|\tau\le t).$$
In particular, the feedbacks may be either reduced to differential form via
the introducing of the intention fields or simply postulated to be differential.
Thus, we constructed an interactive game from the initial ordinary game.
One may use $\varphi(\tau)$ as well as $\varphi^\circ_{[t_0]}(\tau)$ 
in the feedbacks. 

Note that the controls of the first player may be also considered as
interactive if the corrections to the predictions are taken into account
when controls are chosen.
              
The obtained construction may be used in practice to make more adequate
predictions. Namely, {\sl a posteriori\/} analysis of the differential
interactive games allows to make the short-term predictions in such games
[6]. One should use such predictions instead of the initial ones. Note that
at the moment $t_0$ the first player knows the pure controls of other
players at the interval $[t_0,t_0+\Delta t]$ whereas their real freedom
is interpreted as an interactivity of their controls $u_i(t)$. So it is
reasonable to choose $\Delta t$ not greater than the admissible time depth 
of the short-term predictions. Estimations for this depth were proposed
in [6]. 

Na{\"\i}vely, the proposed idea to improve the predictions is to consider
deviations of the real behaviour of players from the predicted ones as
a result of the interactivity, then to make the short-term predictions
taking the interactivity into account and, thus, to receive the 
corrections to the initial predictions. Such corrections may be regarded
as ``psychological'' though really they are a result of different methods
of predictions used by players.

\remark{Remark 5}
The interpretation of the ordinary differential game as an interactive game
also allows to perform the strategical analysis of interactive games. Indeed,
let us consider an arbitrary differential interactive game $A$. Specifying its
$\varepsilon$-representation one is able to construct the associated 
ordinary differential game $B$ with the doubled number of players. Making
some predictions in the game $B$ one transform it back into an interactive 
game $C$. Combination of the strategical long-term predictions in the game 
$B$ with the short-term predictions in $C$ is often sufficient to obtain the 
adequate strategical prognosis for $A$.
\endremark

\remark{Remark 6} The interpretation of the ordinary differential game as
an interactive game is especially useful in situations, when the goals of
players are not known precisely to each other and some more or less rough
suppositions are made.
\endremark

\head II. Dialogues and verbalizable games\endhead

\subhead 2.1. Dialogues\endsubhead
Let us now expose the interactive game formalism for a description of
dialogues as psycholinguistic phenomena [3]. First of all, note that one is
able to consider interactive games of discrete time as well as interactive
games of continuous time above.

\definition{Defintion 3A (the na{\"\i}ve definition of dialogues) [3]}
The {\it dialogue\/} is a 2-person interactive game of discrete time with 
intention fields of continuous time.
\enddefinition

The states and the controls of a dialogue correspond to the speech whereas 
the intention fields describe the understanding. 

Let us give the formal mathematical definition of dialogues now.

\definition{Definition 3B (the formal definition of dialogues) [3]}
The {\it dialogue\/} is a 2-person interactive game of discrete time of 
the form
$$\varphi_n=\Phi(\varphi_{n-1},\vec v_n,\xi(\tau)| 
t_{n-1}\!\leqslant\!\tau\!\leqslant\!t_n).\tag4$$
Here $\varphi_n\!=\!\varphi(t_n)$ are the states of the system at the
moments $t_n$ ($t_0\!<\!t_1\!<\!t_2\!<\!\ldots\!<\!t_n\!<\!\ldots$), 
$\vec v_n\!=\!\vec v(t_n)\!=\!(v_1(t_n),v_2(t_n))$ are the interactive 
controls at the same moments; $\xi(\tau)$ are the intention fields of 
continuous time with evolution equations
$$\dot\xi(t)=\Xi(\xi(t),\vec u(t)),\tag5$$
where $\vec u(t)=(u_1(t),u_2(t))$ are continuous interactive controls with 
$\varepsilon$--represented couplings of feedbacks:
$$u_i(t)=u_i(u_i^\circ(t),\xi(t);\varepsilon_i(t)).$$
The states $\varphi_n$ and the interactive controls $\vec v_n$ are certain
{\sl known\/} functions of the form
$$\aligned
\varphi_n=&\varphi_n(\vec\varepsilon(\tau),\xi(\tau)| 
t_{n-1}\!\leqslant\!\tau\!\leqslant\!t_n),\\
\vec v_n=&\vec v_n(\vec u^\circ(\tau),\xi(\tau)|
t_{n-1}\!\leqslant\!\tau\!\leqslant\!t_n).
\endaligned\tag6
$$
\enddefinition

Note that the most nontrivial part of mathematical formalization of dialogues
is the claim that the states of the dialogue (which describe a speech) are 
certain ``mean values'' of the $\varepsilon$--parameters of the intention
fields (which describe the understanding).

\remark{Important}
The definition of dialogue may be generalized on arbitrary number of players
and below we shall consider any number $n$ of them, e.g. $n=1$ or $n=3$, 
though it slightly contradicts to the common meaning of the word ``dialogue''.
\endremark

\subhead 2.2. The verbalization of interactive games\endsubhead
An embedding of dialogues into the interactive game theoretical picture
generates the reciprocal problem: how to interpret an arbitrary differential
interactive game as a dialogue. Such interpretation will be called the
{\it verbalization}.

\definition{Definition 4 [3]} A differential interactive game of the form
$$\dot\varphi(t)=\Phi(\varphi(t),\vec u(t))$$
with $\varepsilon$--represented couplings of feedbacks 
$$u_i(t)=u_i(u^\circ_i(t),\varphi(t),\dot\varphi(t),\ddot\varphi(t),\ldots
\varphi^{(k)}(t);\varepsilon_i(t))$$
is called {\it verbalizable\/} if there exist {\sl a posteriori\/}
partition $t_0\!<\!t_1\!<\!t_2\!<\!\ldots\!<\!t_n\!<\!\ldots$ and the
integrodifferential functionals
$$\aligned
\omega_n&(\vec\varepsilon(\tau),\varphi(\tau)|
t_{n-1}\!\leqslant\!\tau\!\leqslant\!t_n),\\
\vec v_n&(\vec u^\circ(\tau),\varphi(\tau)|
t_{n-1}\!\leqslant\!\tau\!\leqslant\!t_n)
\endaligned\tag7$$ 
such that
$$\omega_n=\Omega(\omega_{n-1},v_n;\varphi(\tau)|
t_{n-1}\!\leqslant\!\tau\!\leqslant\!t_n).
\tag 8$$
\enddefinition

The verbalizable differential interactive games realize a dialogue in sense
of Def.3.

\remark{Remark 7} One may include $\omega_n$ explicitely into the evolution
equations for $\varphi$
$$\dot\varphi(\tau)=\Phi(\varphi(\tau),\vec u(\tau);\omega_n),\quad 
\tau\in[t_n,t_{n+1}].$$
as well as into the feedbacks and their couplings.
\endremark

The main heuristic hypothesis is that all differential interactive games
``which appear in practice'' are verbalizable. The verbalization means that 
the states of a differential interactive game are interpreted as intention 
fields of a hidden dialogue and the problem is to describe such dialogue 
completely. If a differential interactive game is verbalizable one 
is able to consider many linguistic (e.g. the formal grammar of a related 
hidden dialogue) or psycholinguistic (e.g. the dynamical correlation of 
various implications) aspects of it.

During the verbalization it is a problem to determine the moments $t_i$. A 
way to the solution lies in the structure of $\varepsilon$-representation.
Let the space $E$ of all admissible values of $\varepsilon$-parameters be
a CW-complex. Then $t_i$ are just the moments of transition of the 
$\varepsilon$-parameters to a new cell. 

\subhead 2.3. Psycholinguistic encoding/decoding of dialogues\endsubhead
Let's now describe one practically valuable procedure. If one has a dialogue
it is possible to consider it as a verbalizable game and, hence, as an
interactive game. Then one is able to forget the initial dialogue structure
of this game and to make an attempt to verbalize it. Certainly, a different
dialogue may be obtained in such way. This procedure will be called the {\it
psycholinguistic encoding\/} of the initial dialogue, whereas the reciprocal
one will be called the {\it psycholinguistic decoding}. 

In practice, the psycholinguistic encoding/decoding may include the change
of the communication medium, e.g. the phonic-verbal dialogue may be encoded
as visual-figurative one (the {\it illustration\/}) and vice versa (the
{\it ecphrasis\/}). The analysis of psycholinguistic encoding/decoding may
be useful for a clarification of nature of the speech-script dualism of 
language as well as for an understanding of linguistic aspects of various
forms of art such as dance, martial arts, etc.

\head III. Tactical games\endhead

Tactics as it will be defined below is derived from two independent concepts:
the parametric interactive games and external controls on one hand and
the comments to dialogues on another hand.

\subhead 3.1. Parametric interactive games and external controls\endsubhead
An interactive game may depend on the additional paramaters. Such dependence
is of two forms. First, parameters may appear in the evolution equations:
$$\dot\varphi=\Phi(\varphi,u_1,u_2,\ldots,u_n;\lambda).\tag9A$$
Here, $\lambda$ is a collective notation for parameters. Second, parameters
may appear in feedbacks:
$$u_i(t)=u_i(u_i^\circ(t),\varphi(t),\ldots,\varphi^{(k)}(t);\lambda).\tag9B$$
The dependence of $u_i$ on $\lambda$ is either unknown (incompletely known) or 
known. The least means that $\partial u_i/\partial\lambda$ may be expressed via 
$u_i$ as a function of other variables (such expression are integrodifferential 
on these variables). Both variants of parametric dependence of interactive game 
may be combined together.

The additional paramaters may realize the external controls. In this
situation they depend on time:
$$\lambda=\lambda(t).$$
In practice, such situation appear in the teaching systems. The parameters
are interpreted as controls of a teacher. This example is rather typical.
It shows that the controls $\lambda(t)$ may be considered as ``slow''
whereas the interactive controls $u_i(t)$ as ``quick''.

Of course, one is able to introduce the slow controls $\lambda(t)$, which
belong to the same players as the interactive controls $u_i(t)$ or to their
coalitions. And, certainly, the slow controls of discrete time may be
considered. One may suspect that the discrete time controls $\lambda_n$
realize a convenient approximation for the slow controls $\lambda(t)$, which
is timer in practice.

The slow controls may be interactive.

If dependence of $u_i$ on $\lambda$ is known and one consider the 
$\varepsilon$-representation of feedbacks it is either postulated that 
$\varepsilon$-parameters do not depend on $\lambda$ or claimed that 
$\partial\varepsilon/\partial\lambda$ is expressed via $\varepsilon$ as 
a function of other arguments.
                                                   
\subhead 3.2. Dialogue and comments\endsubhead
Let us consider a $n$-person dialogue
$$\omega_n=\Omega(\omega_{n-1},\vec v_n,\xi(\tau)|
t_{n-1}\!\leqslant\!\tau\!\leqslant\!t_n)$$
with the discrete time interactive controls $\vec v_n$ and the intention
fields governed by the evolution equations
$$\dot\xi(t)=\Xi(\xi(t),\vec u(t)),$$
where $\vec u(t)$ are the continuous interactive controls with
$\varepsilon$--represented couplings of feedbacks:
$$u_i(t)=u_i(u_i^\circ(t),\xi(t);\varepsilon_i(t)).$$
The states $\varphi_n$ and the interactive controls $\vec v_n$ are
expressed as
$$\aligned
\omega_n=&\omega_n(\vec\varepsilon(\tau),\xi(\tau)|
t_{n-1}\!\leqslant\!\tau\!\leqslant\!t_n),\\
\vec v_n=&\vec v_n(\vec u^\circ(\tau),\xi(\tau)|
t_{n-1}\!\leqslant\!\tau\!\leqslant\!t_n).
\endaligned
$$

The discrete time {\it comments\/} $\vartheta_n$ to
the dialogue are defined recurrently as
$$\vartheta_n=\Theta(\vartheta_{n-1},\omega_n,v_n).\tag10$$

Comments to the dialogue at the fixed moment $t_n$ contain various
information on the dialogue. For instance, one may to raise a problem
to restore some features of a dialogue from certain comments or alternatively
what features of a dialogue may be restored from the fixed comment.

The main difference of the comments $\vartheta_n$ from the states $\omega_n$
is the absence of expressions of the first via $\vec\varepsilon(\tau)$ and
$\xi(\tau)$ ($t_{n-1}\!\leqslant\!\tau\!\leqslant\!t_n$).

Comments are applied to the verbalizable games in the same way.
                                                               
\subhead 3.3. Tactical games\endsubhead
Tactical games combine both mechanisms defined above.

\definition{Definition 5} The {\it tactical game\/} is a parametric
verbalizable game with comments, in which the parameters are of discrete
time and coincide with the comments.
\enddefinition

It is really wonderful that such simple definition is applicable to a very
huge class of phenomena. However, it is so! As it was marked above virtually
all known forms of human activity such as scientific researches and
economics, sport or military actions, fine or martial arts, literature,
music and theatre, psychotherapy and even magic may be regarded as certain 
tactical games. Trying to improve the model I have no found any wider concept,
whose using is necessary and effective, whereas the notion of tactical
game may describe these phenomena very correctly.

\remark{Remark 8} The pairs $(v_n,\vartheta_n)$ of discrete time interactive
controls and the comments will be called the {\it tactical actions}, whereas
the continuous time interactive controls $u_i(t)$ will be called the {\it 
instant actions}. The tactical actions may be involved as in the evolution 
equations as in the interactivity.
\endremark

\remark{Remark 9} The description(-construction) of interactive games in 
sense of [4] is a tactical procedure.
\endremark

\subhead 3.4. Applications to computer games (tactical action games and the 
modelling of tactical behavior)\endsubhead
Though some computer games claim that they are tactical it is not so. Really
tactics is not involved in the rules of such games and is accidental.
Otherwords, a player may perform the tactical actions but they are not 
obligatory. The real actions are instant only and their tactical 
interpretation is not substantial. Moreover, artificial players do not perform 
any really tactical actions except the actions following a behavioral pattern. 
The precise definition of a tactical game allows to create, first, the tactical 
action games, where tactical actions are substantial, second, to perform the
computer modelling of tactical behavior for the artificial players. It is not
difficult and is very interesting to do. In present, the author effectively
uses the tactical game constructions in elaboration of tactical computer
games. However, the discussion of practical questions is a bit beyond the
aim of this article.

\head IV. Behavioral self-organization\endhead

This paragraph is devoted to the tactical game modelling of the behavioral
self-organization.

\subhead 4.1. Tactical behavioral self-organization\endsubhead
The tactical properties of players may change during the game. Note that
the tactics is defined by the dependence of the evolution equations and
feedbacks on the comments and by the function $\Theta$, which determines
the comments recurrently. One may fix the dependence on the comments and
vary the function $\Theta$. The procedures of improvement of $\Theta$, which
goal is an adaptation of a player, form the {\it tactical behavioral 
self-organization}. Thus, the tactical behavioral self-organization is a
form of functional self-organization. Numerous methods of modelling of
such self-organization may be, therefore, adopted.

Tactical behavioral self-organization should have a lot of real applications
in all forms of human activity, where tactics appears. There are several
directions of such applications: (1) the teaching, (2) the human adaptation
and self-regulation systems, (3) semi-artificial human-computer systems, 
(4) purely artificial computer systems. The most of them should be based on 
the results, which will be obtained in the computer games.

\subhead 4.2. Applications to computer games (artificial intelligence and
tactical self-organization, tactical RPG)\endsubhead
Combinations of the artificial intelligence and tactical self-organization
are applicable to computer games. In this way the truly tactical RPG may be
created. Note that virtually all known RPG are not tactical and therefore
do not use the tactical behavioral self-organization (in fact, they use
any behavioral self-organization only episodically, the most elaborated forms
of self-organization are applied by the Japanese but they are not widely
distributed). The author, who actively works in this direction, suspects that 
it has a lot of perspectives.

\head V. Tactics and behavioral self-organization for perception
games\endhead

This paragraph is devoted to applications of tactical games to the perception
processes and related subjects (memory, recollection, image understanding,
imagination) as well as to the computer vision and pattern recognition (the
dynamical formation of patterns and perception models during perception as a
result of its self-organization).

\subhead 5.1. Perception games and kaleidoscope-roulettes [7]\endsubhead
Perception processes was understood in the interactive game theoretical terms
in the articles [7]. The main concept is one of the perception game, which
is exposed below. 

\definition{Definition 6} The {\it perception game\/} is a multistage
verbalizable game (no matter finite or infinite) for which the intervals 
$[t_i,t_{i+1}]$ are just the sets. The conditions of their finishing 
depends only on the current value of $\varphi$ and the state of $\omega$ 
at the beginning of the set. The initial position of the set is the final 
position of the preceeding one.
\enddefinition

Practically, the definition describes the discrete character of the
perception and the image understanding. For example, the goal of a concrete
set may be to perceive or to understand certain detail of the whole image.
Another example is a continuous perception of the moving or changing object.

Note that the definition of perception games is applicable to various forms 
of perception, though the most interesting one is the visual perception.
The proposed definition allows to take into account the dialogical character
of the image understanding and to consider the visual perception, image
understanding and the verbal (and nonverbal) dialogues together. It may be
extremely useful for the analysis of collective perception, understanding 
and controlling processes in the dynamical environments -- sports, dancings, 
martial arts, the collective controlling of moving objects, etc. On the other 
hand this definition explicates the self-organizing features of human 
perception, which may be unraveled by the game theoretical analysis. And, 
finally, the definition put a basis for a systematical application of the 
linguistic (e.g. formal grammars) and psycholinguistic methods to the 
image understanding as a verbalizable interactive game with a mathematical 
rigor. 

Let's now consider an interesting class of perception games, the
kaleidoscope-rolettes.

The kaleidoscope-roulette is a result of the attempt to combine the
kaleidoscope, one of the simplest and effective visual game, with the
roulette essentially using the elements of randomness and the treatment of
resonances. The main idea is to substitute random sequences of roulette
by the quasirandom sequences, which may be generated by the interactive
kaleidoscope. The obtained formal definition is below.

\definition{Definition 7} {\it Kaleidoscope-roulette\/} is a perception
game with a quasirandom sequence of quantites $\{\omega_n\}$.
\enddefinition

Certainly, the explicit form of functionals (7) is not known to the players.

Many concrete versions of kaleidoscope-roulettes are constructed. Though 
they are naturally associated with an entertainment their real applications 
may be far beyond it due to their origin and the abstract character of their 
definition. 

Kaleidoscope-roulettes are very interesting due to the resonance
phenomena, which may appear in them. Indeed, though the sequence
$\{\omega_n\}$ is quasirandom the equations (8) for them may have
the resonance solutions. The resonance means a dynamical correlation
of two quasirandom sequences $\{v_n\}$ and $\{\omega_n\}$ whatever $\varphi$
is realized. In such case the quantities $\{v_n\}$ may be comprehended as
``fortune'', what is not senseless in contrast to the ordinary roulette.
However, $v_n$ are {\sl interactive\/} controls and their explicit 
dependence on $\vec u^\circ$ and $\varphi$ is not known. Nevertheless, one 
is able to use {\sl a posteriori\/} analysis and short-term predictions
based on it (cf.[6]) if the time interval $\Delta t$ in the short-term
predictions is not less than the interval $t_{n+1}-t_n$. To do it 
one should slightly improve constructions of [6] to take the discrete-time
character of $v_n$ into account. It allows to perform the short-term 
controlling of the resonances in a kaleidoscope-roulette if they are observed.
The conditions of applicability of short-time predictions to the controlling
of resonances may be expressed in the following form: one should claim that
{\sl variations of the interactivity should be slower than the change of sets 
in the considered multistage game}. 

\remark{Remark 10} The possibility to control resonances by $v_n$ using its
short-term predictions does not contradict to its quasirandomness, because
$v_n$ is quasirandom with respect to $v_{n-1}$ but not to $\varphi(\tau)$
($\tau\in[t_n,t_{n+1}]$).
\endremark

\subhead 5.2. Tactics and behavioral self-organization for perception
games. Perceptive oracles\endsubhead
Many ``representative'' mechanisms in perception processes have the tactical
origin. One should include the memory, the recollection, the image
understanding and the imagination. These phenomena may be described in terms
of the {\it tactical perception games}. Procedures of tactical behavioral
self-organization may be applied to the memory strengthening, the
intensification of creativity during imagination and other psychological
problems.
         
Below we shall discuss an interesting example related to the
kaleidoscope-roulettes, the perceptive oracle.

\definition{Definition 8} The parametric kaleidoscope-roulette, which is
a tactical game, is called a {\it perceptive oracle}. 
\enddefinition

Thus, the perceptive oracle is a kaleidoscope-roulette if the comments are
frozen and in this case the sequence $\{\omega_n\}$ is quasirandom. However,
it may be not so if the comments form tactical actions. In the resonances
the sequence $\{\omega_n\}$ may have some laws. The weaker forms of relations
may also appear, e.g. the integrals $K_\alpha(\omega_n,\vartheta_n,
\omega_{n-1},\vartheta_{n-1},\ldots,\omega_{n-k},\vartheta_{n-k})\equiv 
k_\alpha$ may exist for several subsequent $n$. This explain the choice of
the name "perceptive oracle". The appearing of resonances is manifested
by the omens [5:Rem.7]. So such omens allow to perform some predictions in
the perceptive oracle.

\remark{Remark 11} The perceptive oracle may be regarded as a new tactical
game technique and, thus, may be used for the generation of other tactical
game techniques for the computer games. For instance, it will be very
intersting to combine perceptive oracle with the conveniently generalized 
domino technique.
\endremark

Some words should be said on other applications of the tactical perception 
game formalism. Besides the phenomena of human perception such games may be 
used for elaboration of man-made systems. Thus, the tactical perception games 
are a natural framework for some problems of the computer vision and pattern 
recognition, e.g. the dynamical formation of patterns and perception models 
during perception as a result of its self-organization.

\head VI. Tactical interaction, tactical synthesis and tactical localization 
of control systems\endhead

Note that in a tactical game the comments form a control system with the
initially fixed control algorithm (and tactical behavioral self-organization 
is interpreted as a self-developping of this algorithm). Indeed, if one
consider the control system as unity of structural and functional aspects,
the generation of comments is a structural aspect whereas their representation
as parameters of evolution or interactivity forms a functional aspect.
Application of the tactical game formalism to the simultaneous functioning of 
several control systems (their interaction, synthesis and localization) is 
discussed below.

\subhead 6.1. Tactical interaction of control systems\endsubhead
Let us consider two control systems represented as tactical games defined by
the evolution equations
$$\dot\varphi_1=\Phi_1(\varphi_1,\vec u_1;\vartheta_1)$$
and 
$$\dot\varphi_2=\Phi_2(\varphi_2,\vec u_2;\vartheta_2)$$
with $\varepsilon$--represented couplings of feedbacks
$$u_{1,i}=u_{1,i}(u^\circ_{1,i},\varphi_1,\dot\varphi_1,
\ddot\varphi_1,\ldots \varphi^{(k)}_1;\varepsilon_{1,i},
\vartheta_1)$$
and
$$u_{2,i}=u_{2,i}(u^\circ_{2,i},\varphi_2,\dot\varphi_2,
\ddot\varphi_2,\ldots \varphi^{(k)}_2;\varepsilon_{2,i},
\vartheta_2).$$
The integrodifferential functionals (7) have the form
$$\aligned
\omega_{j,n}&(\vec\varepsilon_j(\tau),\varphi_j(\tau)|
t_{n-1}\!\leqslant\!\tau\!\leqslant\!t_n),\\
\vec v_{j,n}&(\vec u_j^\circ(\tau),\varphi_j(\tau)|
t_{n-1}\!\leqslant\!\tau\!\leqslant\!t_n)
\endaligned$$
and the relations (8)
$$\omega_{j,n}=\Omega_j(\omega_{j,n-1},v_{j,n};\varphi_j(\tau)|
t_{n-1}\!\leqslant\!\tau\!\leqslant\!t_n)$$
hold ($j=1,2$). The comments $\vartheta_1$ and $\vartheta_2$ are defined
recurrently as
$$\vartheta_{1,n}=\Theta_1(\vartheta_{1,n-1},\omega_{1,n},v_{1,n})$$
and
$$\vartheta_{2,n}=\Theta_2(\vartheta_{2,n-1},\omega_{2,n},v_{2,n}).$$

\define\tni{\operatorname{int}}
The {\it tactical interaction\/} is realized by the addition of the
interaction terms into the recurrent formulas for $\vartheta_j$ to produce
the interdetermination of
comments:
$$\vartheta_{1,n}=\Theta_1(\vartheta_{1,n-1},\omega_{1,n},v_{1,n})+
\tilde\Theta_{1,2}^{\tni}(\vartheta_{1,n-1},\vartheta_{2,n-1},\omega_{1,n},
v_{1,n})$$
and
$$\vartheta_{2,n}=\Theta_2(\vartheta_{2,n-1},\omega_{2,n},v_{2,n})+
\tilde\Theta_{2,1}^{\tni}(\vartheta_{2,n-1},\vartheta_{1,n-1},\omega_{2,n},
v_{2,n}).$$
         
\subhead 6.2. Tactical synthesis of control systems\endsubhead
Let us consider $N$ control systems represented as tactical games defined
by the evolution equations
$$\dot\varphi_j=\Phi_j(\varphi_j,\vec u_j;\vartheta_j)$$
($j=1,2,\ldots N$) with $\varepsilon$--represented couplings of feedbacks
$$u_{j,i}=u_{j,i}(u^\circ_{j,i},\varphi_j,\dot\varphi_j,
\ddot\varphi_j,\ldots \varphi^{(k)}_j;\varepsilon_{j,i},
\vartheta_j).$$
The integrodifferential functionals (7) have the form
$$\aligned
\omega_{j,n}&(\vec\varepsilon_j(\tau),\varphi_j(\tau)|
t_{n-1}\!\leqslant\!\tau\!\leqslant\!t_n),\\
\vec v_{j,n}&(\vec u_j^\circ(\tau),\varphi_j(\tau)|
t_{n-1}\!\leqslant\!\tau\!\leqslant\!t_n)
\endaligned$$
and the relations (8)
$$\omega_{j,n}=\Omega_j(\omega_{j,n-1},v_{j,n};\varphi_j(\tau)|
t_{n-1}\!\leqslant\!\tau\!\leqslant\!t_n)$$
hold. The comments $\vartheta_j$ are defined recurrently as
$$\vartheta_{j,n}=\Theta_j(\vartheta_{j,n-1},\omega_{j,n},v_{j,n}).$$

The {\it tactical synthesis\/} is realized by the redefinition of the
recurrent formulas for $\vartheta_j$ to produce the unification of
comments:
$$\vartheta_{j,n}=\tilde\Theta_j(\vartheta_{1,n-1},\ldots\vartheta_{N,n-1},
\omega_{1,n},\ldots,\omega_{N,n},v_{1,n},\ldots v_{N,n}).$$
The functions $\tilde{\boldsymbol\Theta}\!=\!(\tilde\Theta_1,\ldots\tilde
\Theta_N)$ determines the synthesis. It may has various internal structure, 
which is characterized by the set of real arguments of functions 
$\tilde\Theta_j$ and their hierarchical structure. It presupposes that 
$\tilde\Theta_j$ depend not on all triples $(\vartheta_j,\omega_j,v_j)$ and 
various triples may appear in $\tilde\Theta_j$ in coalitions of different form 
and nature. One may think that the functions $\tilde{\boldsymbol\Theta}\!=
\!(\tilde\Theta_1,\ldots\tilde\Theta_N)$ are constructed from the functions
$\boldsymbol\Theta=(\Theta_1,\ldots\Theta_N)$ using some operations, which
realize the synthesis. Thus, the synthesis may be performed in several steps
and be described by an algorithm. Optimization problems for the synthesis
and its construction naturally arise.

\subhead 6.3. Tactical localization of control systems\endsubhead
Tactical localization is a procedure reciprocal to the tactical synthesis.
Let us consider a control system represnted as a tactical game
$$\dot\varphi=\Phi(\varphi,\vec u;\vartheta)$$
with $\varepsilon$--represented couplings of feedbacks
$$u_i=u_i(u^\circ_i,\varphi,\dot\varphi,
\ddot\varphi,\ldots \varphi^{(k)};\varepsilon_i,
\vartheta).$$
The integrodifferential functionals (7) have the form
$$\aligned
\omega_n&(\vec\varepsilon(\tau),\varphi(\tau)|
t_{n-1}\!\leqslant\!\tau\!\leqslant\!t_n),\\
\vec v_n&(\vec u^\circ(\tau),\varphi(\tau)|
t_{n-1}\!\leqslant\!\tau\!\leqslant\!t_n)
\endaligned$$
and the relations (8)
$$\omega_n=\Omega(\omega_{n-1},v_n;\varphi(\tau)|
t_{n-1}\!\leqslant\!\tau\!\leqslant\!t_n)$$
hold. The comments $\vartheta$ are defined recurrently as
$$\vartheta_n=\Theta(\vartheta_{n-1},\omega_n,v_n).$$
To perform the {\it tactical localization\/} means to represent the
control system as a result of a tactical synthesis of $N$ partial control
systems.

The localization of control systems is often used to minimize the resource
expences and to reduce the information circuits. The tactical localization
may be useful in this way.

\head Appendix A. Operative computer games\endhead

The appendix is devoted to the operative computer games and the user
programming of operative units in a multi-user online operative computer
game.

\subhead A.1. Tactics and operative art. Operative computer games\endsubhead
The concept of the operative game has its origin in the Russian tradition of
military science, where operative art is concerned as the intermediate
component of military art between tactics and strategy. Structurally, the
difference between tactics and operative art is that tactics investigates
the controlling of the interactively controlled systems as it was described
in the main text of this article, whereas operative art has deal with the 
controlling of such systems, which are simultaneously the control ones (and 
their controls may be also interactive). Therefore, the operative units are 
in fact the complex tactical formations. 

The most natural models for the operative art are the operative computer 
games (OCG), when computer is modelling the behavior of various operative 
units and is functioning as the control system. The player sends interactively
the commands to the operative units, which are performed by the computer. 
However, a way of their realization is not known completely to the player.
For instance, the distribution of functions between tactical units as well as 
tactical features of these units may be unknown incompletely.

Note that tactics has deal with tactical units whereas the operative art has
deal with operative units, which are tactical formations.

Usually, operative units are represented as schematic non-figurative signs 
for the tactical formations. Thus, the game field is realized in two
different separate views: one has, first, the observable field of figurative 
tactical units and, second, the schematic field of the corresponding 
operative units on the screen of the monitor.

\subhead A.2. Programming of operative units in operative computer 
games\endsubhead
It is an interesting problem of the user programming of operative units in
a multi-user online operative computer game. Ordinarily, such game is
a multi-stage one, so the programming of operative units is done between
the sets. There are two variants. First, the game does not provide the
players by any sources for such programming. In this case an individual
player should write his/her own algorithms for operative units and then
to adapt them to the game. Second, the game may include the unified sources
for the operative unit programming. For instance, it can specify the
programming language and special libraries. Certainly, the language should
be simple, close to the ordinary programming languages, compatible with
them and effective. Undoubtly, none of elaborated programming languages
is not specially adapted for the problem, however, separate features
are useful. The compatibility essentially restricts the effectiveness
because known languages were not created for the multi-user online games. 
However, it will be a problem of the whole computer game community if 
somebody decides to include the programming of operative units in the 
proposed game. I may only hope that the game ``for the programmers'' will 
be interesting enough to compensate such social difficulties and will
stimulate a collective activity on the boundary of computer games and 
the programming art.
\newpage

\Refs
\roster
\item"[1]" Isaaks R., Differential games. Wiley, New York, 1965;\newline
Owen G., Game theory, Saunders, Philadelphia, 1968;\newline
Vorob'ev N.N., Current state of the game theory, Uspekhi Matem. Nauk 15(2) 
(1970) 80-140 [in Russian].
\item"[2]" Juriev D., Interactive games and representation theory. I,II.
E-prints: math.FA/9803020, math.RT/9808098.
\item"[3]" Juriev D., Interactive games, dialogues and the verbalization.
E-print: math.OC/9903001.
\item"[4]" Juriev D., On the description of a class of physical interactive
information systems [in Russian]. Report RCMPI-96/05 [e-version: 
mp\_arc/96-459].
\item"[5]" Juriev D., Games, predictions, interactivity. E-print:
math.OC/9906107.
\item"[6]" Juriev D., The laced interactive games and their {\sl a
posteriori\/} analysis. E-print:\linebreak math.OC/9901043; Differential
interactive games: The short-term predictions. E-print: math.OC/9901074.
\item"[7]" Juriev D., Perceptive games, the image understanding and 
interpretational geometry. E-print: math.OC/9905165; Kaleidoscope-roulette:
the resonance phenomena in perception games. E-print: math.OC/9905180.
\endroster
\endRefs
\enddocument